\documentclass[11pt,leqno]{article}
\usepackage{amsfonts}

\begin{document}


\def\R{{\mathbb R}}
\def\T{{\mathbb T}}
\def\D{{\mathbb D}}
\def\S{{\mathbb S}}
\def\C{{\mathbb C}}
\def\Z{{\mathbb Z}}
\def\N{{\mathbb N}}
\def\H{{\mathbb H}}
\def\B{{\mathbb B}}
\def\diam{\mbox{\rm diam}}
\def\sn{\S^{n-1}}
\def\rr{{\cal R}}
\def\mt{{\Lambda}}
\def\e{\emptyset}
\def\dQ{\partial Q}
\def\dk{\partial K}
\def\endofproof{{\rule{6pt}{6pt}}}
\def\di{\displaystyle}
\def\dist{\mbox{\rm dist}}
\def\sa+{\Sigma_A^+}
\def\du{\frac{\partial}{\partial u}}
\def\dv{\frac{\partial}{\partial v}}
\def\dt{\frac{d}{d t}}
\def\dx{\frac{\partial}{\partial x}}
\def\con{\mbox{\rm const }}
\def\nn{{\cal N}}
\def\mm{{\cal M}}
\def\kk{{\cal K}}
\def\ll{{\cal L}}
\def\vv{{\cal V}}
\def\bb{{\cal B}}
\def\ff{{\cal F}}
\def\i{{\bf i}}
\def\tt{{\cal T}}
\def\uu{{\cal U}}
\def\wloc{W_{\epsilon}}
\def\Int{\mbox{\rm Int}}
\def\dist{\mbox{\rm dist}}
\def\pr{\mbox{\rm pr}}
\def\pp{{\cal P}}
\def\aa{{\cal A}}
\def\cc{{\cal C}}
\def\supp{\mbox{\rm supp}}
\def\Arg{\mbox{\rm Arg}}
\def\In{\mbox{\rm Int}}
\def\con{\mbox{\rm const}\;}
\def\Re{\mbox{\rm Re}}
\def\li{\mbox{\rm li}} 
\def\Seo{S^*_\epsilon(\Omega)}
\def\sdk{S^*_{\dk}(\Omega)}
\def\lae{\Lambda_{\epsilon}}
\def\ep{\epsilon}
\def\oo{{\cal O}}
\def\be{\begin{equation}}
\def\ee{\end{equation}}
\def\beqn{\begin{eqnarray*}}
\def\eeqn{\end{eqnarray*}}
\def\Pr{\mbox{\rm Pr}}

\def\ii{{\imath }}
\def\jj{{\jmath }}
\def\II{{\cal I}}
\def\dd{{\cal D}}
\def\la{\langle}
\def\ra{\rangle}
\def\bs{\bigskip}
\def\Con{\mbox{\rm Const}\;}
\def\do{\partial \Omega}
\def\dk{\partial K}
\def\dl{\partial L}
\def\ll{{\cal L}}
\def\kk{{\cal K}}
\def\kk{{\cal K}}
\def\pr{{\rm pr}}
\def\ff{{\cal F}}
\def\G{{\cal G}}
\def\C{{\bf C}}
\def\dist{{\rm dist}}
\def\dds{\frac{d}{ds}}
\def\con{{\rm const}\;}
\def\Con{{\rm Const}\;}
\def\di{\displaystyle}
\def\oo{\mbox{\rm O}}
\def\hess{\mbox{\rm Hess}}
\def\gi{\gamma^{(i)}}
\def\endofproof{{\rule{6pt}{6pt}}}
\def\ep{\epsilon}
\def\ms{\medskip}
\def\ex{\mbox{\rm extd}}

\def\clip{C^{\mbox{\footnotesize \rm Lip}}}
\def\wlocs{W^s_{\mbox{\footnote\rm loc}}}
\def\Lip{\mbox{\rm Lip}}
\def\lip{\mbox{{\footnotesize\rm Lip}}}
\def\Vol{\mbox{\rm Vol}}
\def\id{\mbox{\rm id}}
\def\co{\; \stackrel{\circ}{C}}

\def\oV{\overline{V}}
\def\saa{\Sigma^+_A}
\def\sa{\Sigma_A}
\def\mta{\Lambda(A, \tau)}
\def\mtaa{\Lambda^+(A, \tau)}

\def\supp{\mbox{\rm supp}}
\def\Arg{\mbox{\rm Arg}}
\def\In{\mbox{\rm Int}}
\def\diam{\mbox{\rm diam}}
\def\e{\emptyset}
\def\endofproof{{\rule{6pt}{6pt}}}
\def\di{\displaystyle}
\def\dist{\mbox{\rm dist}}
\def\con{\mbox{\rm const }}
\def\Box{\spadesuit}
\def\Int{\mbox{\rm Int}}
\def\dist{\mbox{\rm dist}}
\def\pr{\mbox{\rm pr}}
\def\be{\begin{equation}}
\def\ee{\end{equation}}
\def\beqn{\begin{eqnarray*}}
\def\eeqn{\end{eqnarray*}}
\def\la{\langle}
\def\ra{\rangle}
\def\bs{\bigskip}
\def\Con{\mbox{\rm Const}\;}
\def\clip{C^{\mbox{\footnotesize \rm Lip}}}
\def\wlocs{W^s_{\mbox{\footnote\rm loc}}}
\def\Lip{\mbox{\rm Lip}}
\def\lip{\mbox{\footnotesize\rm Lip}}
\def\Re{\mbox{\rm Re}}
\def\li{\mbox{\rm li}} 
\def\ep{\epsilon}
\def\ms{\medskip}
\def\dds{\frac{d}{ds}}
\def\oo{\mbox{\rm O}}
\def\hess{\mbox{\rm Hess}}
\def\id{\mbox{\rm id}}
\def\ii{{\imath }}
\def\jj{{\jmath }}
\def\graph{\mbox{\rm graph}}
\def\span{\mbox{\rm span}}

\def\i{{\bf i}}
\def\C{{\bf C}}

\def\ss{{\cal S}}
\def\tt{{\cal T}}
\def\E{{\cal E}}
\def\rr{{\cal R}}
\def\nn{{\cal N}}
\def\mm{{\cal M}}
\def\kk{{\cal K}}
\def\ll{{\cal L}}
\def\vv{{\cal V}}
\def\ff{{\cal F}}
\def\hh{{\cal H}}
\def\tt{{\cal T}}
\def\uu{{\cal U}}
\def\cc{{\cal C}}
\def\pp{{\cal P}}
\def\aa{{\cal A}}
\def\oo{{\cal O}}
\def\II{{\cal I}}
\def\dd{{\cal D}}
\def\ll{{\cal L}}
\def\ff{{\cal F}}
\def\G{{\cal G}}

\def\hs{\hat{s}}
\def\hz{\hat{z}}
\def\hL{\hat{L}}
\def\hl{\hat{l}}
\def\hl{\hat{l}}
\def\hc{\hat{\cc}}
\def\hbb{\widehat{\cal B}}
\def\hu{\hat{u}}
\def\hX{\hat{X}}
\def\hx{\hat{x}}
\def\hu{\hat{u}}
\def\hv{\hat{v}}
\def\hQ{\hat{Q}}
\def\hC{\widehat{C}}
\def\hF{\hat{F}}
\def\hf{\hat{f}}
\def\hii{\hat{\ii}}
\def\hr{\hat{r}}
\def\hq{\hat{q}}
\def\hy{\hat{y}}
\def\hZ{\widehat{Z}}
\def\hz{\hat{z}}
\def\hE{\widehat{E}}
\def\hR{\widehat{R}}
\def\hell{\hat{\ell}}
\def\hs{\hat{s}}
\def\hW{\widehat{W}}
\def\hS{\widehat{S}}
\def\hV{\widehat{V}}
\def\hB{\widehat{B}}
\def\hhh{\widehat{\cal H}}
\def\hK{\widehat{K}}
\def\hU{\widehat{U}}
\def\hhh{\widehat{\hh}}
\def\hdd{\widehat{\dd}}
\def\hZ{\widehat{Z}}
\def\hal{\hat{\alpha}}
\def\hbe{\hat{\beta}}
\def\hg{\hat{\gamma}}
\def\hrho{\hat{\rho}}
\def\hd{\hat{\delta}}
\def\hphi{\hat{\phi}}
\def\hmu{\hat{\mu}}
\def\hnu{\hat{\nu}}
\def\hsi{\hat{\sigma}}
\def\htau{\hat{\tau}}
\def\hpi{\hat{\pi}}
\def\hep{\hat{\epsilon}}
\def\hxi{\hat{\xi}}
\def\hLa{\widehat{\Lambda}}
\def\hPhi{\widehat{\Phi}}
\def\hPsi{\widehat{\Psi}}

\def\tc{\tilde{C}}
\def\tg{\tilde{\gamma}}  
\def\tV{\widetilde{V}}
\def\tC{\widetilde{\cc}}
\def\tr{\tilde{R}}
\def\tb{\tilde{b}}
\def\tt{\tilde{t}}
\def\tx{\tilde{x}}
\def\tp{\tilde{p}}
\def\tz{\tilde{Z}}
\def\tZ{\tilde{Z}}
\def\tF{\tilde{F}}
\def\tK{\widetilde{K}}
\def\tf{\tilde{f}}
\def\tp{\tilde{p}}
\def\te{\tilde{e}}
\def\tv{\tilde{v}}
\def\tu{\tilde{u}}
\def\tw{\tilde{w}}
\def\ts{\tilde{\sigma}}
\def\tr{\tilde{r}}
\def\tU{\tilde{U}}
\def\tS{\tilde{S}}
\def\tP{\widetilde{\Pi}}
\def\ttau{\tilde{\tau}}
\def\tLip{\widetilde{\Lip}}
\def\tz{\tilde{z}}
\def\tS{\tilde{S}}
\def\tts{\tilde{\sigma}}

\def\tPsi{\tilde{\Psi}}
 \def\tp{\tilde{p}}
\def\tR{\tilde{R}}
\def\tQ{\tilde{Q}}
\def\oL{\tilde{\Lambda}}
\def\tq{\tilde{q}}
\def\tx{\tilde{x}}
\def\ty{\tilde{y}}
\def\tz{\tilde{z}}
\def\tmt{\tilde{\Lambda}}
\def\tg{\tilde{g}}
\def\tsi{\tilde{\sigma}}
\def\ttt{\tilde{t}}
\def\tC{\tilde{C}}
\def\tc{\tilde{c}}
\def\tell{\tilde{\ell}}
\def\trho{\tilde{\rho}}
\def\ts{\tilde{s}}
\def\tB{\widetilde{B}}
\def\thh{\widetilde{\cal H}}
\def\tV{\widetilde{V}}
\def\trr{\tilde{r}}
\def\tv{\tilde{v}}
\def\tu{\tilde{u}}
\def\tw{\tilde{w}}
\def\trho{\tilde{\rho}}
\def\tell{\tilde{\ell}}
\def\tz{\tilde{Z}}
\def\tF{\tilde{F}}
\def\tf{\tilde{f}}
\def\tp{\tilde{p}}
\def\ttau{\tilde{\tau}}
\def\tz{\tilde{z}}
\def\tg{\tilde{\gamma}}  
\def\tV{\widetilde{V}}
\def\tC{\widetilde{\cc}}
\def\tLa{\widetilde{\Lambda}}
\def\tR{\widetilde{R}}
\def\tr{\tilde{r}}
\def\tc{\widetilde{C}}
\def\tD{\widetilde{D}}
\def\tt{\tilde{t}}
\def\tx{\tilde{x}}
\def\tp{\tilde{p}}
\def\tS{\tilde{S}}
\def\tts{\tilde{\sigma}}
\def\tZ{\widetilde{Z}}
\def\tdelta{\tilde{\delta}}
\def\th{\tilde{h}}
\def\tB{\widetilde{B}}
\def\thh{\widetilde{\hh}}
\def\tep{\tilde{\ep}}
\def\tE{\widetilde{E}}
\def\tu{\tilde{u}}
\def\txi{\tilde{\xi}}
\def\teta{\tilde{\eta}}
\def\tnu{\tilde{\nu}}
\def\tphi{\tilde{\phi}}
\def\tOm{\widetilde{\Omega}}
\def\tmt{\widetilde{\mt}}
\def\tPhi{\widetilde{\Phi}}
\def\tal{\tilde{\alpha}}

\def\sr{{\sc r}}
\def\mt{{\Lambda}}
\def\do{\partial \Omega}
\def\dk{\partial K}
\def\dl{\partial L}
\def\wloc{W_{\epsilon}}
\def\span{\mbox{\rm span}}
\def\Jac{\mbox{\rm Jac}}
\def\Vol{\mbox{\rm Vol}}
\def\hh{{\mathcal H}}

\begin{center}
{\Large\bf TANGENT BUNDLES TO REGULAR BASIC SETS IN HYPERBOLIC DYNAMICS} 
\end{center}

\begin{center}
{\sc Luchezar Stoyanov}
\end{center}

\footnotesize

\noindent
{\sc Abstract.} Given a locally maximal compact invariant hyperbolic set $\mt$ for a $C^1$ flow or diffeomorphism
on a Riemann manifold with $C^1$ unstable laminations, we construct an invariant continuous bundle
of tangent vectors to local unstable manifolds that locally  approximates  $\mt$ in a certain way.

\normalsize

\section{\sc Introduction}

Let $M$ be a $C^1$ complete (not necessarily compact)  Riemann manifold $M$ and let $\mt$ be a basic set for 
a $C^1$ flow $\phi_t : M \longrightarrow M$ or a $C^1$ diffeomorphism $f : U \longrightarrow f(U) \subset M$,
where $U$ is an open neighbourhood of $\mt$ in $M$. 
Let $E^u(x)$ and $E^s(x)$ ($x\in \mt$)  be the tangent spaces to the local stable and unstable 
manifolds $W^s_\ep(x)$  and $W^u_\ep(x)$ of size $\ep > 0$, respectively (see Sect. 2).  We say that
$\phi_t$ {\it has $C^1$ unstable laminations over} $\mt$ if the map $\mt \ni x \mapsto E^u(x)$ is $C^1$, i.e. at each $x\in \mt$ 
this map has a linearization that depends continuously on $x$.

Given $x\in \mt$, let $k_x \geq 1$ be the minimal integer such that there exist $\delta_x \in (0,\ep)$ and
a $k_x$-dimensional $C^1$ submanifold $S_x$ of $W^u_{\ep}(x)$ such that 
$\mt \cap W^u_{\delta_x}(x) \subset S_x$. We will then say that $S_x$ is a $u$-submanifold of  
{\it minimal dimension at} $x$ and {\it size} $\delta_x$. 
In general $S_x$ is not (locally) unique, however it turns out that, under certain regularity conditions, 
its tangent space $T_xS_x$ is uniquely determined, invariant  and continuous:


\bs

\noindent
{\bf Theorem 1.} {\it  Assume that $\mt$ is a basic set for the flow $\phi_t$ (or the diffeomorphism $f$) 
and the local  unstable laminations over $\mt$ are  $C^1$. Then there 
exists an integer $k = k^u\geq 1$ and a continuous $d\phi_t$-invariant (or $df$-invariant, respectively)
distribution $E^u_{\mt}(x)$,  $x\in \mt$, of $k$-dimensional linear spaces such that for any 
$x\in \mt$ we have  $k_x = k$, $E^u_\mt(x) \subset E^u(x)$ and $T_xS_x = E^u_\mt(x)$  for any  
$u$-submanifold $S_x$ of minimal dimension at $x$. Moreover, the distribution $E^u_\mt$ is invariant 
under the linearizations of the local stable  holonomy maps, as well.}

\bs

We refer the reader to Sect. 2 for the definition of the local stable holonomy maps.

In a similar way one defines $s$-submanifolds of minimal dimension at $x \in \mt$, and a result similar to the 
above holds for these.

\ms

\noindent
{\bf Remark 1.} It is easy to see that all  vectors in $E^u(x)$ that are tangent to the basic set $\mt$ belong to $E^u_\mt(x)$. 
However, in general the span of these tangent vectors could be a proper subspace  of $E^u_\mt(x)$ -- see Example 2 in
Sect. 4 below.

\ms

The motivation for the introduction of the bundles $E^{u,s}_\mt$ comes from \cite{kn:St2} which deals with Dolgopyat type
spectral estimates (\cite{kn:D}) for Axiom A flows on basic sets, where one has to take into account some fine geometric 
properties of basic sets. As it turns out there, in certain cases when the basic set $\mt$ is `relatively small', considering the
tangent maps $d\phi_t$ over the whole bundle $E^u$ (or $E^s$) does not say enough about the dynamics of the flow
near $\mt$. So, it appears,  restricting $d\phi_t$ over the smallest possible invariant subbbundle of $E^u$ 
would be beneficial\footnote{Although in \cite{kn:St2} we found a way around the problem by different means.}, and it is natural
to believe that this might be the case in some other situations, as well, particularly when the dimension of the spaces
$E^u(x)$ (or $E^s(x)$) is large.

In Sect. 3 below we prove Theorem 1 for flows; for diffeomorphisms the argument is very similar. Some
basic definitions are given in Sect. 2, while Sect. 4 contains two examples.

\section{\sc Preliminaries}

Let $M$ be a $C^1$ complete (not necessarily compact) 
Riemann manifold,  and $\phi_t : M \longrightarrow M$ ($t\in \R$) a $C^1$ flow on $M$. A
$\phi_t$-invariant closed subset $\mt$ of $M$ is called {\it hyperbolic} if $\mt$ contains
no fixed points  and there exist  constants $C > 0$ and $0 < \lambda < 1$ and a $d\phi_t$-invariant decomposition 
$T_xM = E^0(x) \oplus E^u(x) \oplus E^s(x)$ of $T_xM$ ($x \in \mt$) into a direct sum of non-zero linear subspaces,
where $E^0(x)$ is the one-dimensional subspace determined by the direction of the flow
at $x$, $\| d\phi_t(u)\| \leq C\, \lambda^t\, \|u\|$ for all  $u\in E^s(x)$ and $t\geq 0$, and
$\| d\phi_t(u)\| \leq C\, \lambda^{-t}\, \|u\|$ for all $u\in E^u(x)$ and  $t\leq 0$.
Here $\|\cdot \|$ is the {\it norm} on $T_xM$ determined by the Riemann metric on $M$.

A non-empty compact $\phi_t$-invariant hyperbolic subset $\mt$ of $M$ which is not a single 
closed orbit is called a {\it basic set} for $\phi_t$ if $\phi_t$ is transitive on $\mt$ 
and $\mt$ is locally maximal, i.e. there exists an open neighbourhood $V$ of
$\mt$ in $M$ such that $\mt = \cap_{t\in \R} \phi_t(V)$. 

From now on we will assume that $\mt$ is a basic set for the flow $\phi_t$.
For $x\in \Lambda$ and a sufficiently small $\epsilon > 0$ let 
$$\wloc^s(x) = \{ y\in M : d (\phi_t(x),\phi_t(y)) \leq \epsilon \: \mbox{\rm for all }
\: t \geq 0 \; , \: d (\phi_t(x),\phi_t(y)) \to_{t\to \infty} 0\: \}\; ,$$
$$\wloc^u(x) = \{ y\in M : d (\phi_t(x),\phi_t(y)) \leq \epsilon \: \mbox{\rm for all }
\: t \leq 0 \; , \: d (\phi_t(x),\phi_t(y)) \to_{t\to -\infty} 0\: \}$$
be the (strong) {\it stable} and {\it unstable manifolds} of size $\epsilon$. Then
$E^u(x) = T_x \wloc^u(x)$ and $E^s(x) = T_x \wloc^s(x)$. 


Assuming $\ep_1 \in (0,\ep_0]$  is sufficiently small, for any $x\in \mt$ and 
$y \in \mt\cap W^s_{\ep_1}(x)$  the local {\it holonomy map}
$\hh_x^y : \mt \cap W^u_{\ep_1}(x) \longrightarrow \mt \cap W^u_{\ep_0}(y)$
along stable laminations is well-defined and uniformly H\"older continuous (see e.g.  \cite{kn:PSW}). 
Recall that the map $\hh_x^y$ is defined as follows. Given $z\in \mt \cap W^u_{\ep_1}(x)$, there exist unique
$z' \in W^s_{\ep_1}(z)$ and $y'\in W^u_{\ep_0}(y)$  such that $\phi_t(z') = y'$ for some $t\in \R$, $|t|\leq\ep_0$.
Then we set $\hh_x^y(z) = y'$.  Under the additional condition that the unstable laminations are $C^1$, the maps $\hh_x^y$ are
$C^1$ as well (see e.g. Fact (2) on  p. 647 in \cite{kn:Ha}). That is, for each $z\in \mt\cap W^u_{\ep_1}(x)$
the map $\hh_x^y$ has a {\it linearization} $L_x^y(z) : E^u(z) \longrightarrow E^u(\hh_x^y(z))$ 
at $z \in \mt \cap W^u_{\ep_1}(x)$ and $\|L_x^y(z)\| \leq C$ for some constant $C > 0$ 
independent of $x$, $y$ and $z$. 

\section{\sc Proof of Theorem 1}


So, assume that $\phi_t : M \longrightarrow M$ ($t\in \R$) is a $C^1$ flow on $M$ and 
$\mt$ is a basic set for $\phi_t$.

By Whitney's Theorem (see e.g. \cite{kn:M}), the stable holonomy map $\hh_x^y$ has a $C^1$ extension 
$\thh_x^y : W^u_{\ep_1}(x) \longrightarrow W^u_{\ep_0}(y)$ 
with $\|d \thh_x^y\| \leq C$ on $W^u_{\ep_1}(x)$. Since 
$\|(d \thh_x^y(x))^{-1}\| = \|(d \hh_x^y(x))^{-1}\| = \|d \hh_y^x(y)\| \leq C$,
it follows from the Inverse Function Theorem that there exists $\ep_2 \in (0,\ep_1]$ such that 
$$\thh_x^y : W^u_{\ep_2}(x) \longrightarrow  \thh_x^y (W^u_{\ep_2}(x)) \subset W^u_{\ep_0}(y)$$
is a diffeomorphism for all $x\in \mt$ and all $y\in \mt \cap W^s_{\ep_1}(x)$. 

Let $k_x$ ($x\in \mt$) be the integers introduced in Sect. 1. Set
$k = \min_{x\in \mt} k_x$,
and fix an arbitrary $x_0 \in \mt$ with $k_{x_0} = k$.

\bs

\noindent
{\bf Lemma  1.} {\it $k_x = k$ for every $x\in \mt$. Moreover, there exist
$\delta > 0$ and for every $x\in \mt$ a $k$-dimensional $C^1$ submanifold $S_x$
of $W^u_{\delta}(x)$ such that $\mt \cap W^u_{\delta}(x) \subset S_x$.}

\bs

\noindent
{\it Proof.} First, notice that if $S_x$ is a  $C^1$ submanifold 
of $W^u_{\delta}(x)$ of minimal dimension such that $\mt \cap W^u_{\delta}(x) \subset S_x$ 
for some  $x\in \mt$, then for any $y \in \mt\cap W^s_{\ep_1}(x)$, $\thh_x^y(S_x)$
is a $C^1$ submanifold of $W^u_{\ep_0}(y)$ with 
$\mt \cap W^u_{\delta/C}(y) \subset \thh_x^y(S_x)$, so $k_y \leq k_x$.
Similarly, $k_x\leq k_y$, so $k_x = k_y$. Moreover, we can choose the same $\delta > 0$
for all $y \in \mt \cap W^s_{\ep_1}(x)$. 

Let  $S = S_{x_0}$ be a  $u$-submanifold  of $W^u_{\delta}(x_0)$ of minimal dimension 
such that $\mt \cap W^u_{\delta}(x_0) \subset S$ for some $\delta > 0$. Given
$z\in \mt$, for every sufficiently large $t > 0$ we have  
$\phi_t(\mt \cap W^u_{\delta/2}(x_0)) \cap W^s_{\ep_1}(z) \neq \e$ (see e.g. \cite{kn:KH}).
Take $t > 0$ with this property so large that $\phi_t$ is expanding on $W^u_{\delta/2}(x_0)$, i.e.
$d(\phi_t(p), \phi_t(q)) \geq d(x,y)$ for all $p,q \in W^u_{\delta/2}(x_0)$.
Let $z' = \phi_t(x) \in W^s_{\ep_1}(z)$ for some $x \in \mt \cap W^u_{\delta/2}(x_0)$. 
Since $S_x = S\cap W^u_{\delta/2}(x) \subset S\cap W^u_{\delta}(x_0)$, $S_x$ contains
$\mt \cap W^u_{\delta/2}(x)$ and $\dim(S_x) = k \leq k_x$, it follows that $k_x = k$ and
$S_x$ is a $u$-submanifold of minimal dimension at $x$ and size $\delta$.

It now follows that $k_{z'} = k$ and, since $\phi_t$ is expanding on $W^u_{\delta/2}(x)$, 
there exists a $u$-submanifold of minimal dimension at $z'$ of size $\delta/2$. The remark
in the beginning of the proof now shows that $k_z = k_{z'} = k$, and
there exists a $u$-submanifold of minimal dimension at $z$ of size $\delta/(2C)$.
\endofproof

\bs

\noindent
{\bf Lemma 2.} {\it For any $x\in \mt$ and any two $u$-submanifolds $S$ and $S'$ of minimal dimension at $x$
we have $T_x S = T_x S'$. }

\bs

\noindent
{\it Proof.} Let $S$ and $S'$ be $u$-submanifolds  of minimal dimension $k$ at $x$
for some $x\in \mt$. Take $\delta > 0$ so small that $\mt\cap W^u_{\delta}(x) \subset S$ and
$\mt\cap W^u_{\delta}(x) \subset S'$.

Assume that $T_xS \setminus T_xS' \neq \e$ and fix an arbitrary
$v \in T_xS \setminus T_xS'$. Then we must have $n = \dim(W_{\delta}^u(x)) > k$. 
Choosing a smaller $\delta > 0$ if necessary and using an appropriate
submanifold chart for $S'$ about $x$ in $W^u_{\delta}(x)$ one constructs a
$C^1$ submanifold $S''$ of $W^u_{\delta}(x)$ such that $\dim(S'') = n-1$,
$S' \subset S''$ and $v \notin T_xS''$. Then $T_xS + T_xS'' = E^u(x) = T_x(W^u_{\delta}(x))$, so
the submanifolds $S$ and $S''$ are transversal at $x$ (see e.g. \cite{kn:GP}). Consequently,
there exists $\delta' \in (0,\delta]$ so that $S''' = S\cap S'' \cap W^u_{\delta'}(x)$ is a submanifold
of $W^u_{\delta'}(x)$ of dimension at most $k-1$. Since 
$\mt \cap W^u_{\delta'}(x) \subset S'''$, this is a contradiction with the minimality of $k = k_x$.

Thus we must have $T_xS \subset T_xS' $ and similarly $T_xS' \subset T_xS$.
\endofproof

\bs

\noindent
{\it Proof of Theorem 1.} For any $x\in \mt$ define $E^u_\mt(x) = T_xS$, where $S$ is an arbitrary 
$u$-submanifold of minimal dimension at $x$. By Lemma 2, the definition is correct. Moreover,
$E^u_\mt(x) \subset E^u(x)$ and $T_xS_x = E^u_\mt(x)$  for any  $u$-submanifold $S_x$ of  
minimal dimension at $x$. The latter implies $E^u_\mt (\phi_t(x)) = d\phi_t(x) \cdot E^u_\mt(x)$
for any $t\in \R$ and any $x\in \mt$.

It remains to show that the map $\Phi: \: \mt\ni x \mapsto E^u_\mt(x)$ is continuous. 
Fix for a moment $z\in \mt$ and let $S_z$ be an  $u$-submanifold  of  minimal dimension at $z$.
Then for any $t\in \R$ with a sufficiently small $|t|$, $\phi_t(S_z)$ is an $u$-submanifold  of  
minimal dimension at $z'= \phi_t(z)$,
so $E^u_\mt(z') = T_{z'}(\phi_t(S_z)) = d\phi_t(z) \cdot (T_zS_z)$ depends continuously on  
$z' = \phi_t(z)$. That is, $\Phi$ is continuous along the orbit of $z$.

Next, choosing a sufficiently small $\delta > 0$, for $z'\in W^u_\delta(z)\cap \mt$ we have
$z'\in S_z$, so $E^u_\mt(z') = T_{z'}S_z$, and therefore $\Phi(z') = E^u_\mt(z')$ depends continuously on 
$z'\in W^u_\delta(z)\cap \mt$.

Finally, for $z'\in W^s_{\ep_1}(z)\cap \mt$, $\thh_{z}^{z'}(S_z)$ is an $u$-submanifold  of  
minimal dimension at $z'$, so
$E^u_\mt(z') = d\thh_z^{z'}(z)\cdot (E^u_\mt(z))$, and therefore $\Phi(z') = E^u_\mt (z')$ depends continuously
on $z'\in W^s_{\ep_1}(z)\cap \mt$. 

The above and the continuity of the local product near $z$ implies that $\Phi$ is continuous at $z$.
\endofproof

\section{\sc Examples}

We use two very well known kind of flows to give: 1) a non-trivial example where $k^u < \dim (E^u)$, and
2) an example where the span of the tangent vectors in $E^u(z)$ to $\mt$ at $z$ form a 
proper subspace of $E^u_\mt(z)$ (for any $z\in \mt$). We use open billiard flows for the first
example and geodesic flows on hyperbolic manifolds for the second. In higher dimensions both kind of
flows should provide a great variety of non-wandering sets with complicated geometric (and metric)
structures (see e.g \cite{kn:Ka}), however it appears they would not be easy to investigate. 
It is quite possible that one could use other models that are easier to describe to get similar examples. 

There are various possible ways to define tangent vectors to a subset of a manifold. Here we use
the most straightforward one.
Given $z \in \mt$, let $\exp^u_z : E^u(z) \longrightarrow W^u_{\ep_0}(z)$  and
$\exp^s_z : E^s(z) \longrightarrow W^s_{\ep_0}(z)$ be the corresponding
{\it exponential maps}. A  vector $b\in E^u(z)\setminus \{ 0\}$ is called  {\it tangent to $\mt$} at
$z$ if there exist infinite sequences $\{ v^{(m)}\} \subset  E^u(z)$ and $\{ t_m\}\subset \R\setminus \{0\}$
such that $\exp^u_z(t_m\, v^{(m)}) \in \mt \cap W^u_{\ep}(z)$ for all $m$, $v^{(m)} \to b$ and 
$t_m \to 0$ as $m \to \infty$.  It is easy to see that a vector $b\in E^u(z)\setminus \{ 0\}$ is  tangent to $\mt$ at
$z$ iff there exists a $C^1$ curve $z(t)$, $0\leq t \leq a$, in $W^u_{\ep}(z)$ for some $a > 0$ 
with $z(0) = z,\: \dot{z}(0) = b$, and $z(t) \in \mt$ for arbitrarily small $t >0$. Denote by $\hE^u(z)$
the set of all non-zero tangent vectors $v \in E^u(z)$ to $\mt$ at $z$.
Tangent vectors to $\mt$ in $E^s(z)$ and the set $\hE^s(z)$ are defined similarly.  
Clearly the sets $\hE^{u,s}(z)$ are non-empty and $d\phi_t$-invariant and $\hE^{u,s}(z) \subset E^{u,s}(z)$. 
However in general the span of $\hE^{u,s}(z)$ does not coincide with $E^{u,s}(z)$ -- see Example 2 below.

\ms

\noindent
{\bf Open Problem.} Under the assumptions in Sect. 1, assume that $k_u < \dim(E^u(z))$ (or $k_s < \dim(E^s(z))$), $z\in \mt$.
Does there always exist a $C^1$ submanifold $M'$ of $M$ of positive codimension and an open
neighbourhood $U$ of $\mt$ in $M$ such that $\phi_t(U\cap M') \subset M'$ for all $t \in \R$?

\ms

One would expect the answer to be negative, however at this stage we do not have examples to prove this. 
Example 1 below is non-trivial in a certain sense, and in fact might be good enough to answer the above
question in the negative, however we do not go that far here.

\subsection{\sc Open billiard flows}

Let $K$ be a subset of ${\R}^{n}$ ($n\geq 2$) of the form
$K = K_1 \cup K_2 \cup \ldots \cup K_{k_0}$, where $K_i$ are compact
 strictly convex disjoint domains in $\R^{n}$ with 
$C^r$ {\it boundaries} $\Gamma_i = \dk_i$ ($r \geq 2$) and $k_0 \geq 3$. Set 
$\Omega = \overline{{\R}^n \setminus K}$ and $\Gamma = \partial K$.
We assume that $K$ satisfies the following {\it no-eclipse condition}: 
 {\it for every pair $K_i$, $K_j$ of different connected components 
of $K$ the convex hull of $K_i\cup K_j$ has no
common points with any other connected component of $K$}. 
With this condition, the {\it billiard flow} $\phi_t$ defined on the {\it cosphere bundle} $S^*(\Omega)$ 
in the standard way is called an {\it open billiard flow}.
It has singularities, however its restriction to the {\it non-wandering set} $\Lambda$ (the set of those
$x\in S^*(\Omega)$ such that the trajectory $\{ \phi_t(x) : t\in \R\}$ is bounded) has only 
simple discontinuities at reflection points\footnote{Notice that the natural projection of $\phi_t$ on the quotient space
$S^*(\Omega)/\sim$, where $\sim$ is the equivalence relation $(q,v) \sim (p,w)$ iff $q=p$ and $v = w$ or 
$q = p\in \dk$ and $v$ and $w$ are symmetric with respect to $T_q(\dk)$, is continuous. Moreover
whenever both $x$ and $\phi_t(x)$ are in the interior of $S^*(\Omega)$ and sufficiently
close to $\Lambda$, the map $y \mapsto \phi_t(y)$ is smooth on a neighbourhood of $x$.
It follows from well-known results of Sinai that $\Lambda$ is a hyperbolic
set for $\phi_t$, and it is easily seen that $\Lambda$ is the maximal compact $\phi_t$-invariant subset 
of $S^*(\Omega)$. Moreover, it follows from the natural symbolic coding for the natural section of the flow that the 
periodic points are dense in $\Lambda$, and $\phi_t$ is transitive on $\Lambda$.
Thus, $\Lambda$ is a basic set for $\phi_t$ and the classical theory of hyperbolic flows
applies  (see e.g. Part 4 in \cite{kn:KH}).}.  
Moreover, $\Lambda$  is compact and $\phi_t$ is hyperbolic and transitive on $\Lambda$.
Finally, it follows from \cite{kn:St3} that when the minimal distance between  distinct connected 
components of $K$ is relatively large compared to the maximal sectional curvature of $\dk$ the open billiard
flow on $\mt$ satisfies a certain pinching condition which implies that the (un)stable laminations over $\mt$ are $C^1$.

In the following example we will use the natural symbolic coding of the open billiard.
Let $A$ be the $k_0\times k_0$ matrix with entries $A(i,j) = 1$ if $i \neq j$ and $A(i,i)= 0$ for all $i$, and let
$\sa$ be the set of all sequences  $\eta = (\eta_j)_{j=-\infty}^\infty$
of integer numbers  $1\leq \eta_j \leq k_0$ such that  $\eta_{j} \neq \eta_{j+1}$ for all $j \in \Z$.
Let $\pr_1 : S^*(\Omega) \longrightarrow \Omega$ be the natural {\it projection}. 
Given $\xi \in \sa$, let $(P_j(\xi))_{j=-\infty}^\infty$
be the successive reflection points of the unique billiard trajectory in the exterior of $K$ such that
$P_j(\xi) \in K_{\xi_j}$ for all $j \in \Z$ (see e.g. \cite{kn:St1}).  Define the map $\Phi : \sa \longrightarrow \mt \cap S^*_{\dk}(\Omega)$
by $\Phi(\xi) = (P_0(\xi), (P_1(\xi) - P_0(\xi))/ \|P_1(\xi) - P_0(\xi)\|)$. Then $\Phi$ is a bijection
such that $\Phi\circ \sigma = B \circ \Phi$, where $B :  \mt \cap S^*_{\dk}(\Omega) \longrightarrow  \mt \cap S^*_{\dk}(\Omega)$ 
is the {\it billiard ball map} from boundary to boundary, and $\sigma$ is the shift map on $\sa$.

\bs

\noindent
{\bf Example 1.} Assume that $n = 3$ and there exists a plane $\alpha$ 
such that each of the domains $K_j$ is symmetric with respect to $\alpha$.  Setting $K' = K\cap \alpha$ and
$\Omega' = \Omega \cap \alpha$, it is easy to observe that every billiard trajectory generated by
a point in $\mt$ is entirely contained in $\alpha$. That is, $\mt = \mt'$, where $\mt'$ is the 
non-wandering set for the open billiard flow in $\Omega'$. Thus, $\dim(E^u_\mt(z)) = 1 < \dim(E^u(z)) = 2$ 
for any $z \in \mt$. This example is of course trivial, since $\mt$ is contained in the flow-invariant
submanifold $S^*(\Omega')$ of $S^*(\Omega)$.

However with a small local perturbation of the boundary $\dk$ of $K$ we can get a non-trivial example.
Choosing standard cartesian coordinates $x,y,z$ in $\R^3$, we may assume that
$\alpha$ is given by the equation $z = 0$, i.e. $\alpha = \R^2\times \{0\}$.
Let $\pr_1: S^*(\R^3) \sim \R^3 \times \S^2 \longrightarrow \R^3$ be the natural projection, and let
$C = \pr_1(\mt)$. We may choose the coordinates $x,y$ in the plane $\alpha = \{ z=0\}$ so that
the line $y =0$ is tangent to $K'_1$ and $K'_2$ and $K'$ is contained in
the half-plane $y \geq 0$. Let $q_1\in K'_1$ and $q_2\in K'_2$ be such that $[q_1,q_2]$ is the
shortest segment connecting $K'_1$ and $K'_2$. Take a point $q_1' \in \dk'_1$ close to
$q_1$ and such that the $y$-coordinate of $y'_1$ is smaller than that of $q_1$. Consider the {\it open
arc} $\aa$ on $\dk'_1$ connecting $q_1$ and $q_1'$. It is clear that $\aa \cap C = \e$.

Let $f: \R^3\longrightarrow \R^3$ be a $C^1$ (we can make it even $C^\infty$) diffeomorphism with 
$f(x) = x$ for all $x$ outside a small open set $U$ such that $q_1 \in \overline{U}$ 
and $U \cap \dk' \subset \aa$. Then for any $q\in C$ the tangent 
planes $T_q(\dk)$ and  $T_q(\partial \tK)$ coincide. We can choose $f$ so that $\tK_i = f(K_i) = K_i$ 
for $i > 1$,  $\tK_1 = f(K_1)$ is strictly convex, and  $\tnu(f(q)) \notin \alpha$ for $q \in \aa$ arbitrarily close to $q_1$. Here 
$\tnu$ is the outward  unit normal field to $\partial \tK$. 

Notice that the non-wandering set $\tmt$ for the billiard
flow $\tphi_t$ in the closure $\tOm$ of the exterior of $\tK$ in $\R^3$ coincides with $\mt$. Indeed, let
$\tPhi : \sa \longrightarrow \mt \cap S^*_{\partial \tK}(\tOm)$ be the coding map for the billiard
trajectories in $\tOm$ and let $\tz = (\tq,\txi) \in S^*_{\partial \tK}(\tOm)$ belong to $\tmt$. 
Then $\tz = \tPhi(\xi)$ for some $\xi \in \sa$, so $z = \Phi(\xi) \in \mt\cap S^*_{\dk} (\Omega)$. 
Since any reflection point $P_i(\xi) \in C$, by the choice of $f$ we have
that $P_j(\xi)$ are the successive reflection points of a billiard trajectory in $\tOm$ and this 
must be the trajectory determined by $\tPhi(\xi)$. Thus, $\tz = z$, so $\tz \in \mt$. This argument 
also shows that $\mt \subset \tmt$, so $\tmt = \mt \subset S^*(\alpha)$. Thus, 
$\dim(E^u_{\tmt}(z)) = 1 < \dim (E^u(z))$ for any $z \in \tmt$.
However, it is clear from the construction that $S^*(\alpha \cap \Omega)$ is not invariant 
with respect to the billiard flow $\tphi_t$. Moreover, it is not difficult to see that there is 
no two-dimensional submanifold $\tal$ of $\tOm$ such that $S^*(\tal)$ is $d\tphi_t$-invariant and
$\mt \subset S^*(\tal)$ . 

Indeed, assume such $\tal$ exists; then $\tal$ is a union of billiard trajectories and
$C = \pr_1(\mt) \subset \tal$. We will first show that
$\tal = \alpha$ outside a large disk  $D$ containing $K'$. Changing the coordinates $x,y$ in the plane
$\alpha = \{ z=0\}$, we may assume that the line $y =0$ is tangent to $K'_2$ and $K'_3$ and 
$K'$ is contained in the half-plane $y \geq 0$. Let $p_2\in K'_2$ and $p_3\in K'_3$ be such that 
$[p_2,p_3]$ is the shortest segment connecting $K'_2$ and $K'_3$. Clearly $p_2, p_3 \in \tal$.
Since $p_2$ is not an isolated point in $C$, there are points in $C \cap \dk'_2$ arbitrarily close
to $p_2$. All of them are in $\tal$, so the curve $\tal \cap \dk_2$ is tangent to $\dk'_2$ at $p_2$.
Similarly, the curve $\tal \cap \dk_3$ is tangent to $\dk'_3$ at $p_3$.
Take point $p_2' \in \tal \cap \dk_2$ close to $p_2$ and $p_3' \in \tal \cap \dk_3$ close to $p_3$
such that the $y$-coordinates of $p'_2$and $p'_3$ are smaller than those of $p_2$ and $p_3$, respectively. 
Consider the arc $C_2$ on $\tal \cap \dk_2$ connecting $p_2$ and $p_2'$ and the arc $C_3$ on 
$\tal \cap \dk_3$ connecting $p_3$ and $p_3'$. We have $(p_2,v_2) \in \mt$, where
$v_2 = (p_3-p_2)/\|p_3-p_2\|$. For any $p\in C_2$ sufficiently close to $p_2$ and any 
$v \in S^*_p(\tal)$ sufficiently close to $v_2$, the billiard trajectory $\gamma(p,v)$ in $S^*(\tOm)$ 
issued from $(p,v)$ is contained in $S^*(\tal)$, so its projection in $\R^3$ is contained in $\tal$. 
Thus, the first reflection point of $\gamma(p,v)$ belongs to $C_3$. That is, for any $p\in C_2$ 
sufficiently close to $p_2$ and any  $v \in S^*_p(\tal)$ sufficiently close to $v_2$, the straight line
segment ray issued from $p$ in the direction of $v$ intersects the curve $C_3$. 
Similarly, for any $q\in C_3$ sufficiently close to $p_3$ and any  $w \in S^*_q(\tal)$ sufficiently close 
to $-v_2$, the straight line segment ray issued from $q$ in the direction of $w$ intersects the curve 
$C_2$. Replacing the curves $C_2$ and $C_3$ by shorter ones (i.e. replacing the points $p_2'$ and
$p'_3$ by points on $C_2$ and $C_3$ closer to $p_2$ and $p_3$, respectively), we now have that
for any $p\in C_2$ and any $q\in C_3$, the straight line segment $[p,q]$ lies in $\tal$. This shows that
$C_2$ and $C_3$ lie in the same plane, and it is clear from the choice of $C_2$ and $C_3$ that this plane
is tangent to $\alpha$ at $p_2$, so it must coincide with $\alpha$. Moreover, the union of the segments $[p,q]$ 
with $p\in C_2$ and any $q\in C_3$ is part of both $\alpha$ and $\tal$. 

From the latter one easily derives that if $D$ is an open disk in $\alpha$ containing 
$K'$, then $\alpha\setminus D \subset \tal$. Since $\mt$ is a nowhere dense subset of 
$S^*(\Omega)$ (and $S^*(\tOm)$, as well), every point $q\in \alpha\setminus K$ can be
approximated arbitrarily well by points of the form $\pr_1(\phi_t(x,\xi))$, 
where $t > 0$ and $(x,\xi) \in S^*(\alpha\setminus D)$. This implies $\alpha\setminus K \subset \tal$,
and therefore $\alpha\setminus K = \tal$, which is a contradiction  with the definition of
$\tK$ and the perturbation $f$. Hence there does not exist a two-dimensional submanifold $\tal$ of 
$\tOm$ such that $S^*(\tal)$ is $d\tphi_t$-invariant and $\mt \subset S^*(\tal)$ .

\subsection{\sc Geodesic flows on manifolds of negative curvature}

Let $X$ be a complete (not necessarily compact) connected Riemann manifold of constant curvature 
$K = -1$ and dimension $\dim(X) = n+1$, $n \geq 1$,  and let $\phi_t : M = S^*(X) \longrightarrow M$ 
be the geodesic flow on the {\it unit cosphere bundle} of $X$.  According to a classical result of 
Killing and Hopf
any such $X$ is a {\it hyperbolic manifold}, i.e. $X$ is isometric to $\H^{n+1}/\Gamma$, where 
$$\H^{n+1} = \{ (x_1, \ldots, x_{1}) \in \R^{n+1} : x_{1} > 0\}$$
is the upper half-space in $\R^{n+1}$ with the {\it Poincar\'e metric} 
$ds^2(x) = \frac{1}{x^2_{1}} (dx_1^2 + \ldots + dx_{1}^2)$
and $\Gamma$ is a {\it Kleinian group}, i.e. a discrete group of isometries (M\"obius transformations) 
acting freely and discontinuously on $\H^{n+1}$. See e.g. \cite{kn:Ratc}  for basic information on hyperbolic manifolds.  
Given a hyperbolic manifold $X = \H^{n+1}/\Gamma$,
the {\it limit set} $L(\Gamma)$ is defined as the set of accumulation points of all $\Gamma$ orbits in 
$\overline{\partial \H^{n+1}}$, the topological closure of $\partial \H^{n+1} = \{ 0\} \times \R^n$ 
including $\infty$.

Throughout this sub-section we will assume that $\Gamma$ is torsion-free and finitely generated (then
$\Gamma$ is called {\it geometrically finite}) and {\it non-elementary}, i.e. $L(\Gamma)$ is infinite
(then $L(\Gamma)$ is a closed non-empty nowhere dense subset of 
$\partial \overline{\H^{n+1}}$ without isolated points; see e.g. Sect. 12.1 in \cite{kn:Ratc}).
A geometrically finite Kleinian group with no parabolic elements is called 
{\it convex cocompact}. If $X$ is compact, then $\Gamma$ is called a {\it cocompact lattice}.

The {\it non-wandering set}  $\mt$ of  $\varphi_t : M  \longrightarrow M$ 
(also known as the {\it convex core} of $X = \H^{n+1}/\Gamma$)
is the image in $M$ of the set of all points of $S^*(\H^{n+1})$ generating geodesics with end 
points in $L(\Gamma)$. When $\Gamma$ is convex cocompact, the non-wandering set $\mt$  is compact.

Notice that the class of hyperbolic manifolds $X = \H^{n+1}/\Gamma$, with $\Gamma$ a non-elementary convex cocompact Kleinian group, contains all classical and non-classical Schottky manifolds (cf. e.g. Sect. 12.1 in \cite{kn:Ratc}). In this case
the (un)stable laminations of the geodesic flow over $\mt$ are  always $C^1$ (in fact $C^\infty$).

\bs

\noindent
{\bf Example 2.} Consider the case $n = 2$ and let $\H^{3}/\Gamma$ be a hyperbolic manifold generated by 
a non-elementary convex cocompact Kleinian group $\Gamma$ with 
an "Apolonian packing" limit set $L(\Gamma)$ (see e.g. \cite{kn:Su}). Then
$L(\Gamma)$ is a subset of $\R^2$ such that for every $x\in L(\Gamma)$ the set of tangent vectors to
$L(\Gamma)$ at $x$ is one-dimensional, while $L(\Gamma) \cap U$ is not contained in an one-dimensional
submanifold of $\R^2$ for any open neighbourhood $U$ of $x$ in $\R^2$. 
Given $z\in \mt$, the local unstable manifold $W^u_\ep(z)$ is given by the projection into $\H^{3}/\Gamma$ of
the outward normal field to a (part of) a horosphere $S$ in $\H^3$. The latter are just spheres in $\H^3$ tangent to $\partial \H^3$
or planes in the interior of $\H^3$ parallel to $\partial \H^3$. Any $z'\in \mt \cap W^u_\ep(z)$ can be identified
with the point $x' \in L(\Gamma)$ which is just the limit at $+\infty$ of the geodesic determined by $z'$. 
It then follows that  the set of tangent vectors in $E^u(z)$ to $\mt$ at $z$ is one-dimensional, while $\mt \cap U$ is not contained 
in an one-dimensional submanifold of $W^u_\ep(z)$ for any $\ep > 0$, i.e. $k^u = 2$.

\bs

\bs

{\sc University of Western Australia, Crawley WA 6009, Australia}

{\it E-mail address:}  stoyanov@maths.uwa.edu.au

\end{document}